%% Amstex file, this is the version of April 9-10, 2016
%% Revision Belgrade June 3, 2016; June 6-11, 2016.
%% New, enlarged version starting June 13-16, 18-19, June, 2016, Belgrade.
%% New version on June 21, 23, 29, 30, 2016. Version August 7, 2016, Krf.
%% Corrections August 16, 2016, Belgrade. Further corrections 21.08.2016.
%% Theorem 1 changed on August 22-23, 2016. Small revision August 16, 2016.
%% Belgrade, September 16, 2016. F. Gora, 17.09.2016. Belgrade, Sept. 22, 2016,
%% Sept. 25, October 1-2, 2016.

\magnification1200
\input amstex.tex
\documentstyle{amsppt}
\nopagenumbers
\hsize=12.5cm
\vsize=18cm
\hoffset=1cm
\voffset=2cm

\footline={\hss{\vbox to 2cm{\vfil\hbox{\rm\folio}}}\hss}

\def\DJ{\leavevmode\setbox0=\hbox{D}\kern0pt\rlap
{\kern.04em\raise.188\ht0\hbox{-}}D}

\baselineskip=13pt
\def\hf{{\textstyle{1\over2}}}
\def\a{\alpha}\def\b{\beta}
\def\d{{\,\roman d}}
\def\e{\varepsilon}\def\E{{\roman e}}

\def\b{\beta} \def\g{\gamma}

\def\t{\theta}
\def\={\;=\;}

\def\zt{\zeta(\hf+it)}

\def\E{{\roman e}}

\def\z{\zeta}

\def\t{\theta}
\def\hf{{\textstyle{1\over2}}}

\def\le{\leqslant} \def\ge{\geqslant}
%%%%%%%%%%% Fonts macros %%%%%%%%%%%%
\font\tenmsb=msbm10
\font\sevenmsb=msbm7
\font\fivemsb=msbm5
\newfam\msbfam
\textfont\msbfam=\tenmsb
\scriptfont\msbfam=\sevenmsb
\scriptscriptfont\msbfam=\fivemsb
\def\Bbb#1{{\fam\msbfam #1}}

\def \NN {\Bbb N}

\font\ff=cmr8

\baselineskip=13pt

\font\teneufm=eufm10
\font\seveneufm=eufm7
\font\fiveeufm=eufm5
\newfam\eufmfam
\textfont\eufmfam=\teneufm
\scriptfont\eufmfam=\seveneufm
\scriptscriptfont\eufmfam=\fiveeufm
\def\mathfrak#1{{\fam\eufmfam\relax#1}}

\font\tenmsb=msbm10
\font\sevenmsb=msbm7
\font\fivemsb=msbm5
\newfam\msbfam
     \textfont\msbfam=\tenmsb
      \scriptfont\msbfam=\sevenmsb
      \scriptscriptfont\msbfam=\fivemsb
\def\Bbb#1{{\fam\msbfam #1}}

\def \NN {\Bbb N}

  \def\rightheadline{{\hfil{\ff
  Some remarks on the differences between ordinates of consecutive zeta zeros}\hfil\tenrm\folio}}

  \def\leftheadline{{\tenrm\folio\hfil{\ff
   Aleksandar Ivi\'c }\hfil}}
  \def\emptyheadline{\hfil}
  \headline{\ifnum\pageno=1 \emptyheadline\else
  \ifodd\pageno \rightheadline \else \leftheadline\fi\fi}

\topmatter
%\nopagenumbers
\title
Some remarks on the differences between ordinates of consecutive zeta zeros
\endtitle
\author   Aleksandar Ivi\'c
 \endauthor

\nopagenumbers

\medskip

\address
Aleksandar Ivi\'c, Katedra Matematike RGF-a
Universiteta u Beogradu, \DJ u\v sina 7, 11000 Beograd, Serbia
\endaddress

\bigskip
\email {
\tt
aleksandar.ivic\@rgf.bg.ac.rs, aivic\_2000\@yahoo.com }\endemail
\dedicatory
\enddedicatory
\abstract
{If $0 < \g_1 \le \g_2 \le \g_3 \le \ldots$ denote ordinates of complex zeros
of the Riemann zeta-function, then several results involving the maximal
order of $\g_{n+1}-\g_n$ and the sum
$$
\sum_{0<\g_n\le T}{(\g_{n+1}-\g_n)}^k  \qquad(k>0)
$$
are proved. }
\endabstract
\endtopmatter

\document
\head
1. Introduction
\endhead
Let $0 < \g_1 \le \g_2 \le \g_3 \le \ldots$ denote ordinates of complex zeros
of the Riemann zeta-function
$$
\z(s) = \sum_{n=1}^\infty n^{-s}\qquad(\Re s > 1).
$$
For $\Re s \le 1$ one defines $\z(s)$ by analytic continuation (see the monographs
of A. Ivi\'c \cite{11} and  E.C. Titchmarsh \cite{19} 
for the properties of $\z(s)$). Here the Riemann Hypothesis (RH), that all complex
zeros of $\z(s)$ satisfy $\Re s = \hf$, is not assumed. Thus if equality among the
$\g_n$'s occurs for some $n$, it does not necessarily mean that the zero
$\rho_n = \b_n + i\g_n$  is not simple,
i.e., $\z(\rho_n) = 0$ and  $\z'(\rho_n) =0$. Namely one could have $\g_n = \g_{n+1}$,
$\rho_n = \b_n + i\g_n, \rho_{n+1} = \b_{n+1} + i\g_{n+1}$ with $\b_n\ne \b_{n+1}$,
and both $\rho_n$ and $\rho_{n+1}$ simple.
Although all numerical evidence points
to the simplicity of all zeta zeros, proving this is an open and difficult question.
In fact, it seems that the simplicity of zeta-zeros and the RH are independent statements
in the sense that, as far as it is known, both statements could be true or false,
or one true and the other one false.

Problems involving $\g_{n+1}-\g_n$, the difference between consecutive ordinates
of the zeros (if the zeros are arranged according to the size of their imaginary parts)
are of great interest. Since $\z(\b_n-i\g_n)=0$ if  $\z(\b_n+i\g_n)=0$, one may consider
without loss of generality that $\g_n>0$ for all $n$.  One of the 
natural problems is to investigate the sum
$$
S_k(T) := \sum_{0<\g_n\le T}{(\g_{n+1}-\g_n)}^k, \leqno(1.1)
$$
where $k$ is a fixed positive number.
 A. Fujii
\cite{5}, \cite{6} proved in 1975 that, for a fixed $k\in\NN$,
$$
C_1\frac{N(T)}{(\log T)^k} \;\le\; S_k(T) \;\le\; C_2\frac{N(T)}{(\log T)^k}.
\leqno(1.2)
$$
In (1.2) we have $0 < C_1 = C_1(k) < C_2 = C_2(k), T \ge T_0 = T_0(k)$, and $N(T)$ denotes
the number of $\g_n$'s not exceeding $T$, counted with multiplicities. Recall that by the classical
 Riemann--von Mangoldt formula (see e.g., 
 \cite{11} or \cite{19}) we have
$$
N(T) = {T\over2\pi}\log\left({T\over2\pi}\right) - {T\over2\pi} + {7\over8}
+ S(T) + O\left({1\over T}\right),\quad S(T) = {1\over\pi}\arg\z(\hf
+ iT).\leqno(1.3)
$$
Here $\arg\z({1\over2}+ iT)$ is obtained by continuous
variation along the segments
joining the points $2, 2 + iT, {1\over2} + iT$, starting with the value 0.
If $T$ is the ordinate of a zero lying on the critical line,
 then $S(T) =S(T + 0)$. One has (see \cite{19}) the bounds
$$
S(T) \ll \log T,\quad S(T) = o(\log T) \quad({\roman {LH}}),\quad
S(T) \ll {\log T\over\log\log T}\quad ({\roman {RH}}),\leqno(1.4)
$$
where LH denotes the (hitherto unproved) Lindel\"of hypothesis that
$$
\zt \;\ll_\e\; |t|^\e.
$$
The LH is a consequence of the RH (see \cite{11} or \cite{19}), but it is not known
whether the converse is true.
Here $f(x) \ll_\e g(x)$ (same as $f(x) = O_\e(g(x))$) means that the implied $\ll$ (or $O$) constant depends
only on $\e$.
The bounds in (1.2) are explicit, but they are stated to hold only if $k\in\NN$. The last restriction
can be easily removed. Indeed, we shall show in Section 3 that (1.2) holds for any fixed $k>1$.

Note that from (1.3) and the first bound in (1.4) we have unconditionally
$$
\g_{n+1}-\g_n \;\ll\; 1.\leqno(1.5)
$$
From Theorem 9.12 of \cite{18} it follows that (1.5) can be improved to
$$
\g_{n+1}-\g_n \;\le\; \frac{A}{\log\log\log\g_n}\qquad(A>0,\, n\ge n_0).\leqno(1.6)
$$
R.R. Hall and W.K. Hayman \cite{10} showed that any constant $A>\pi/2$ is permissible
in (1.6).
Also from (1.3) and (1.4),
 on the RH, the bound (1.6) can be improved to
$$
\g_{n+1}-\g_n \;\ll\; \frac{1}{\log\log\g_n}.\leqno(1.7)
$$

\medskip
The purpose of this article is to investigate $S_k(T)$, as well as some problems involving the
order of $\g_{n+1}-\g_n$ and the frequency of values of $n$ for which this difference is large.
Our results are primarily explicit.

\medskip
{\bf Acknowledgement}. I wish to thank D.R. Heath-Brown and T. Trudgian for valuable remarks.

\medskip
\head
2. The maximal order of $\g_{n+1}-\g_n$
\endhead

\medskip
Although  improving the upper bounds (1.6)  and (1.7) seems difficult, one can derive explicit
bounds, namely replace the $\ll$-constant in (1.7) by an explicit value. This is contained in

\medskip
THEOREM 1.
{\it  Under the Riemann hypothesis one has}
$$
\g_{n+1}-\g_n \;\le\;\bigl(\frac\pi2 +o(1)\bigr)\frac{1}{\log\log\g_n}\qquad(n\to\infty).\leqno(2.1)
$$

\medskip
{\bf Proof}. To prove (2.1), we shall use the  bound, valid under the RH,
$$
|S(T)| \;\le\; \left({1\over4} + o(1)\right){\log T\over\log\log T}
\qquad(T\to\infty).\leqno(2.2)
$$
This is Theorem 2 of  E. Carneiro, V. Chandee and M. Milinovich \cite{3}.
It improves the previous bound of D.A. Goldston and S.M. Gonek \cite{8},
who had the constant $\hf+o(1)$ in (2.2), which yields (2.1) (see their
Corollary 1) with the constant $\pi+ o(1)$.  Actually, in \cite{3} Carneiro et al.
have shown that
$$
|S(T)| \;\le\; \left({1\over4} + O\Bigl(\frac{\log\log\log T}{\log\log T}\Bigr)\right)
{\log T\over\log\log T}.
$$

 We use (1.3) with $T=\g_n, H>0,$ and (2.2) to obtain
$$
\eqalign{
N(T+H) - N(T)&= \frac{1}{2\pi}\int_T^{T+H}\log\left({t\over2\pi}\right)\d t
+ S(T+H)-S(T)+ O\left({1\over T}\right)\cr&
\ge \frac{H}{2\pi}\log\left({T\over2\pi}\right) -
\left({1\over2} + o(1)\right){\log T\over\log\log T}> 0\cr}
$$
for
$$
H \;=\; \Bigl(\frac\pi2 +o(1)\Bigr)\frac{1}{\log\log T}\qquad(T\to\infty).
$$
Thus
$$
\g_{n+1} \in \left[\g_n, \g_n+ \frac{\hf\pi +o(1)}{\log\log\g_n}\,\right],
$$
which implies then (2.1). Clearly the term $o(1)$ in (2.1) can be replaced by 
the more explicit
$$
O\left(\frac{\log\log\log\g_n}{\log\log\g_n}\right).
$$

 \medskip
 We remark that, although the unconditional bound
 (1.5) is weaker than (1.6), one can obtain relatively
 simply an explicit bound for $\g_{n+1}-\g_n$. Namely we take advantage of the recent
 bound of T. Trudgian \cite{20}
$$
|S(T)| \le 0.112\log T + 0.278\log\log T + 2.510\qquad(T \ge \E).\leqno(2.3)
$$
 If we take $T=\g_n$ in (1.3) and use (2.3) we obtain, for $H>0, T\ge T_0$ and 
 some number $\t$ for which $|\t|\le1$,
$$
\eqalign{
N(T+H) - N(T)&= \frac{1}{2\pi}\int_T^{T+H}\log\left({t\over2\pi}\right)\d t
+ S(T+H)-S(T)+ O\left({1\over T}\right)\cr&
\ge \frac{H}{2\pi}\log\left({T\over2\pi}\right) + \t0.225\log T > 0\cr}
$$
for $H = 1.414$ and $T\ge T_0$. This gives then unconditionally
$$
\g_{n+1}-\g_n \;\le\;1.414\qquad(n\ge n_0),\leqno(2.4)
$$
and with some effort one could determine $n_0$ in (2.4) explicitly.

\medskip
Determining the maximal order of $\g_{n+1}-\g_n$ is a difficult problem. Note that from (1.5) one
has unconditionally
$$
\sum_{0<\g_n\le T}(\g_{n+1}-\g_n) = \sum_{0<\g_n\le T, \g_{n+1}\ne\g_n}(\g_{n+1}-\g_n) = T+O(1).
\leqno(2.5)
$$
Hence
$$
T+O(1) = \sum_{0<\g_n\le T}(\g_{n+1}-\g_n) \le
N(T)\max_{0<\g_n\le T}(\g_{n+1}-\g_n),
$$
and from (1.3) one obtains
$$
\max_{0<\g_n\le T}(\g_{n+1}-\g_n) \ge \frac{2\pi\bigl(1+o(1)\bigr)}{\log(T/2\pi)}\qquad
(T\to\infty).\leqno(2.6)
$$
\medskip
The lower bound in (2.6) is quite explicit, but it is weak and probably far from the true order of the
quantity on the left-hand side. In his paper \cite{16}, A.M. Odlyzko  states that
under the GUE (Gaussian Unitary Ensemble hypothesis, see \cite{16} and \cite{12}) it is
plausible that
$$
 \max_{0<\g_n\le T}(\g_{n+1}-\g_n)
\;\sim\; \frac{8}{\sqrt{2\log T}}\quad(T\to\infty).\leqno(2.7)
$$
On the other hand, D. Joyner in \cite{14} brings forth that
under the so-called Dyson--Montgomery hypothesis, explained in  \cite{14}, one has
$$
\max_{0<\g_n\le T}(\g_{n+1}-\g_n) \;\ll\; \frac{1}{\sqrt{\log T\log\log T}}.\leqno(2.8)
$$
Note that (2.7) and (2.8) cannot both be true, since they contradict one another.
The very slow variation of $\sqrt{\log\log T}$ makes a numerical comparison of (2.7)
and (2.8) difficult.

\medskip
\head
3. Some remarks on the moments of $\g_{n+1}-\g_n$
\endhead

\medskip
In this section we shall show  that (1.2) holds for any fixed $k>1$, not
necessarily an integer. We
assume that $k>1$ is fixed and start from (2.5).
Then, by H\"older's inequality,
$$
T + O(1) \le \left\{\sum_{0<\g_n\le T}(\g_{n+1}-\g_n)^{k}\right\}^{1/k}\Bigl\{N(T)\Bigr\}^{1-1/k}.
$$
Since $T \sim 2\pi N(T)/\log T$ by (1.3), the above inequality yields immediately
$$
\sum_{0<\g_n\le T}(\g_{n+1}-\g_n)^k \ge \frac{\bigl(2\pi+o(1)\bigr)^k}{(\log T)^k}N(T)\qquad(T\to\infty).
\leqno(3.1)
$$
Note that (3.1) is the lower bound inequality in (1.2), with the explicit value
$$
C_1 = C_1(k) = \bigl(2\pi+\e\bigr)^k
$$
for any given $\e>0$.

\medskip
To obtain the upper bound inequality, recall that the upper bound in (1.2)
 holds for $k\in \NN$ and suppose that $\a$
satisfies $k < \a < k+1$ for some $k\in\NN$. Then write
$$
\sum_{0<\g_n\le T}(\g_{n+1}-\g_n)^\a = F_1(T) + F_2(T),\leqno(3.2)
$$
say. We have, on using the upper bound in (1.2),
$$
\eqalign{
F_1(T) &:= \sum_{0<\g_n\le T, \g_{n+1}-\g_n \le 1/\log T}(\g_{n+1}-\g_n)^\a\cr&
\le (\log T)^{k-\a}\sum_{0<\g_n\le T}(\g_{n+1}-\g_n)^k \le C_2(k)\frac{N(T)}{(\log T)^\a}.\cr}
$$
Similarly, using (1.2) with $k+1$ in place of $k$, we have
$$
\eqalign{
F_2(T) &:= \sum_{0<\g_n\le T, \g_{n+1}-\g_n > 1/\log T}(\g_{n+1}-\g_n)^\a\cr&
\le (\log T)^{k+1-\a}\sum_{0<\g_n\le T}(\g_{n+1}-\g_n)^{k+1} \le C_2(k+1)\frac{N(T)}{(\log T)^{\a}}.\cr}
$$
Inserting the bounds for $F_1(T)$ and $F_2(T)$ in (3.2) we obtain the desired upper bound for $S_\a(T)$.

\medskip
An asymptotic formula for $S_k(T)$, when $k\ne0,1$, is hard to obtain. One can obtain such a formula
if one assumes the RH
and the Gaussian Unitary Ensemble (GUE) conjecture (see \cite{16} for a detailed account). 
This says that, for
$$
0 \le \a < \b < \infty,\quad \delta_n = \frac{1}{2\pi}(\g_{n+1}-\g_n)\log(\frac{\g_n}{2\pi}),
$$
we have
$$
\sum_{\g_n\le T,\delta_n\in[\a,\b]}1 = \left(\int_\a^\b p(0,u)\d u + o(1)\right)\frac{T}{2\pi}
\log(\frac{T}{2\pi})\qquad(T\to\infty).
$$
 Then one has, as shown by the author in \cite{12},
$$
\sum_{\g_n\le T}(\g_{n+1}-\g_n)^k = \bigl(c_1(k) + o(1)\bigr)
\left(\frac{2\pi}{\log(\frac{T}{2\pi})-1}\right)^{k-1}T\quad(T\to\infty)\leqno(3.3)
$$
for fixed $k\ge0$, thus not necessarily an integer. Here $c_1(0) = c_1(1) = 1$, and in general
$$
c_1(k) := \int_0^\infty p(0,u)u^k\d u,
$$
where $p(0,u)$ is the function appearing in the GUE conjecture.
We have
$$
\eqalign{&
1 - \left(\frac{\sin\pi u}{\pi u}\right)^2 = \sum_{k=0}^\infty p(k,u),
\cr&
p(0,u) = \frac{1}{3}\pi^3 u^2 - \frac{2}{15}\pi^4u^4 +
\frac{1}{315}\pi^6u^6+\cdots\quad(u \to 0+),
\cr&
\log p(0,u) = -\frac{\pi^2}{8} + o(1)\qquad(u\to\infty).
\cr}
$$

From (1.3) one infers that the {\it average} distance $\g_{n+1}-\g_n$ is $2\pi/\log(\g_n/(2\pi))$.
Thus a natural question is to investigate the quantities
$$
\mu := \liminf_{n\to\infty}\frac{\g_{n+1}-\g_n}{2\pi/\log(\g_n/(2\pi))},\quad
\lambda := \limsup_{n\to\infty}\frac{\g_{n+1}-\g_n}{2\pi/\log(\g_n/(2\pi))}.\leqno(3.4)
$$
A Selberg \cite{18} in 1946 indicated (without proof) that $\mu<1$ and $\lambda>1$ holds unconditionally,
but no particular values of $\mu$ and $\lambda$ have been found yet. On the RH, several authors
worked on this problem over the years and produced explicit values of $\mu$ and $\lambda$.
 For example,  Feng and Wu \cite{4}
 obtained the values $\mu \le 0.514$ and $\lambda \ge 2.7327$.
J. Bredberg \cite{1} proved that for sufficiently large $T$ there a subinterval of
$[T, 2T]$ of length at least $2.766\times\frac{2\pi}{\log(T/2\pi)}$ in which $\zt$ does not vanish.
Thus, on the RH, one has $\lambda\ge 2.766$.

\medskip
A stronger variant of (3.4) is that there exist constants $\mu<1$ and $\lambda >1$   such that
$$
\frac{\g_{n+1}-\g_n}{2\pi/\log(\g_n/(2\pi))} \le \mu,\qquad
\frac{\g_{n+1}-\g_n}{2\pi/\log(\g_n/(2\pi))}  \ge \lambda\leqno(3.5)
$$
for a positive proportion of $n$'s. This was stated by A. Fujii in \cite{6}, 
and a detailed proof of (3.5) may be
found on pp. 246-249 of E.C. Titchmarsh's monograph \cite{19}.

\medskip
It is interesting to investigate what is the number of $\g_n$'s not exceeding $T$ for which the
distance $\g_{n+1}-\g_n$ is larger or smaller than the average distance. This problem, and
some related questions, will be discussed in the next section.

\medskip
%\vfill
%\break

\head
4. Lower bounds for sums of large differences of $\g_{n+1}-\g_n$
\endhead
\medskip
We begin our discussion on the frequency of occurrences of $\g_{n+1}-\g_n$.
First note that, for a given positive constant $C$,
$$
\eqalign{
&\sum_{0<\g_n\le T}(\g_{n+1}-\g_n)^2 = \sum_{0<\g_n\le T;\g_{n+1}-\g_n
\le C/\log (T/2\pi) }(\g_{n+1}-\g_n)^2
\cr&\qquad\qquad\qquad\qquad\quad+ \sum_{0<\g_n\le T;\g_{n+1}-\g_n>C/\log (T/2\pi) }
(\g_{n+1}-\g_n)^2\cr&
\le (C^2+o(1))\frac{N(T)}{\log^2T} +
\left(\sum_{0<\g_n\le T;\g_{n+1}-\g_n>\frac{C}{\log(T/2\pi)} }1\right)^{1\over2}
\left(\sum_{\g_n\le T}(\g_{n+1}-\g_n)^4\right)^{1\over2}.
\cr}
$$
Thus it follows, on using (1.2), that
$$
(C_1(2)-C^2+o(1))\frac{N(T)}{\log^2 T} \le
 \left(\sum_{0<\g_n\le T;\g_{n+1}-\g_n>\frac{C}{\log(T/2\pi)} }1\right)^{1/2}
\left(\frac{C_2(4)N(T)}{\log^4 T}\right)^{1/2},
$$
which yields unconditionally
$$
\sum_{0<\g_n\le T;\g_{n+1}-\g_n>C/\log(T/2\pi)}1 \;\ge\; \frac{(C_1(2)-C^2+o(1))^2}{C_2(4)}N(T)
\quad(T\to\infty),\leqno(4.1)
$$
and the bound (4.1) is non-trivial if $0 < C < \sqrt{C_1(2)}$.

\medskip
If one assumes the RH, then the $\g_n$'s are exactly 
the zeros of Hardy's function (see the author's monograph
\cite{13} for an extensive account)
$$
Z(t) := \zt\bigl(\chi(\hf+it)\bigr)^{-1/2}, \;\z(s) = \chi(s)\z(1-s),
$$
which is real-valued and satisfies $|Z(t)| = |\zt|$. Hardy's function is thus an invaluable
tool for the investigation of zeros of $\z(s)$ on the {\it critical line} $\Re s = \hf$.
If one also assumes that almost all the $\g_n$'s are simple,
then (4.1) can be used for obtaining an alternative
proof of Theorem 2 in the paper of Gonek--Ivi\'c
\cite{9}. Following Fujii's arguments one can find numerical values of
the constants in (4.1), but they will certainly produce poor values of the constant in Theorem 2
in \cite{9}. However, the proof of this result assumes both the RH and the Pair Correlation
conjecture, and both of these are strong assumptions.

 The quantity $\log(T/2\pi)$ appearing in (4.1) is natural, because of (1.3) we
already noted that the average spacing between the $\g_n$'s is $2\pi/\log(\g_n/2\pi)$.
Moreover, with increasing $C$ the sum
in (4.1) decreases, so one has to have an expression such as $C_1(2)-C^2$ on the right-hand side of (4.1).

\medskip  In view of (1.3) one can rewrite (3.3) as
$$
\sum_{\g_n\le T}(\g_{n+1}-\g_n)^k = \bigl(c_1(k) + o(1)\bigr)
\left(\frac{2\pi}{\log(\frac{T}{2\pi})-1}\right)^{k}N(T)\quad(T\to\infty).\leqno(4.2)
$$
With this notation (4.1) becomes then
$$
\sum_{0<\g_n\le T;\g_{n+1}-\g_n>C/\log(T/2\pi)}1 \;\ge\;
\frac{\Bigl\{(2\pi)^2c_1(2)-C^2+o(1)\Bigr\}^2}{(2\pi)^2c_1(4)}N(T),\leqno(4.3)
$$
and one has then only to calculate explicitly the values of $c_1(2)$ and $c_1(4)$ and
insert them in (4.3). This will produce an explicit bound in the range
$$
0 < C < 2\pi\sqrt{c_1(2)}.$$

\medskip
A variant of the approach leading to (4.3) is as follows. Recall that we have (2.4), namely
$$
\sum_{0<\g_n\le T}(\g_{n+1}-\g_n) = T+O(1).\leqno(4.4)
$$
Write, for a given $C>0$,
$$
\eqalign{
\sum_{\g_n\le T}(\g_{n+1}-\g_n) &=
\sum_{\g_n\le T;\g_{n+1}-\g_n\le C/\log(T/2\pi)}(\g_{n+1}-\g_n)\cr&
+ \sum_{\g_n\le T;\g_{n+1}-\g_n> C/\log(T/2\pi)}(\g_{n+1}-\g_n)\cr&
= S_1(T;C) + S_2(T;C),\cr}\leqno(4.5)
$$
say. One has trivially
$$
S_1(T;C) \le \frac{C}{\log(T/2\pi)}\sum_{\g_n\le T}1 = \frac{C}{\log(T/2\pi)}N(T).\leqno(4.6)
$$
On the other hand, by the Cauchy-Schwarz inequality, we obtain
$$
S_2(T;C) \le {\left\{\sum_{\g_n\le T;\g_{n+1}-\g_n> \frac{C}{\log(T/2\pi)}}1
\sum_{\g_n\le T;\g_{n+1}-\g_n> \frac{C}{\log(T/2\pi)}}(\g_{n+1}-\g_n)^2 \right\}}^{\frac12}.
\leqno(4.7)
$$
We have trivially
$$
\sum_{\g_n\le T;\g_{n+1}-\g_n> \frac{C}{\log(T/2\pi)}}(\g_{n+1}-\g_n)^2
\le \sum_{\g_n\le T}(\g_{n+1}-\g_n)^2,
$$
and one can estimate the last sum by (1.2). However, if one assumes the Riemann hypothesis, 
then A. Fujii \cite{7} showed that one has
$$
\sum_{\g_n\le T}(\g_{n+1}-\g_n)^2 \le 9\cdot \frac{2\pi T}{\log(T/2\pi)}\quad(T\ge T_0).\leqno(4.8)
$$
Consequently from (4.4)--(4.8), on the RH, we have
$$
\sum_{0<\g_n\le T;\g_{n+1}-\g_n>C/\log(T/2\pi)}1 \ge \frac{T\log(T/2\pi)}{18\pi}\left(1-\frac{C}{2\pi}
\right)^2 + O(T).\leqno(4.9)
$$
Note that (4.9) is an explicit inequality, and it is non-trivial for $0<C< 2\pi$, that is, for the difference
between consecutive ordinates which is smaller than the average difference.

\medskip
In the above two approaches we have exploited the sum in (1.1) with $k=1$ and $k=2$. One can
work with general $k$ in (1.2), but it is unclear which $k$ will yield the best lower bound for the
sum in (4.9).

\medskip
We summarize the preceding discussion in

\medskip
THEOREM 2. {\it With the notation introduced above we have unconditionally, if $0 < C < \sqrt{C_1(2)}$,}
$$
\sum_{0<\g_n\le T;\g_{n+1}-\g_n>C/\log(T/2\pi)}1 \;\ge\;
\frac{(C_1(2)-C^2+o(1))^2}{C_2(4)}N(T)\quad(T\to\infty).
$$
{\it Moreover, if the RH is assumed, then for for $0<C< 2\pi, T\ge T_0>0$ we have}
$$
\sum_{0<\g_n\le T;\g_{n+1}-\g_n>C/\log(T/2\pi)}1 \ge \frac{T\log(T/2\pi)}{18\pi}\left(1-\frac{C}{2\pi}
\right)^2 + O(T).
$$

\medskip
\head
5. Upper  bounds for sums of large differences of $\g_{n+1}-\g_n$
\endhead
A natural problem is to consider upper bounds
 the sum in Theorem 2. An explicit upper bound for this sum is easily obtained.
 Namely, by using (2.5), we have
$$
\eqalign{
\sum_{0<\g_n\le T;\g_{n+1}-\g_n>C/\log(T/2\pi)}1 &\le
\frac{1}{C}\log(T/2\pi)\sum_{0<\g_n\le T}(\g_{n+1}-\g_n) \cr&
= \frac{1}{C}\log(T/2\pi)\Bigl(T+O(1)\Bigr)\cr &
= \frac{2\pi T}{2\pi C}\log(T/2\pi) + O(\log T) \cr&= \frac{2\pi}{C}N(T) + O\Bigl(\frac{T}{C}\Bigr).\cr}
$$
This gives, unconditionally and uniformly for any $C>0$,
$$
\sum_{0<\g_n\le T;\g_{n+1}-\g_n>C/\log(T/2\pi)}1 \le \frac{2\pi}{C}N(T) +
O\Bigl(\frac{T}{C}\Bigr).\leqno(5.1)
$$
 Using the upper bound in (1.2) with general $k$ one obtains similarly
$$
\sum_{0<\g_n\le T;\g_{n+1}-\g_n>C/\log(T/2\pi)}1 \;\le\;\frac{C_2(k)}{C^k}\Bigl(1+o(1)\Bigr)N(T)
\quad(T\to\infty),\leqno(5.2)
$$
but it is unclear for what range of $C$ and the value of $k$ this  bound is optimal.
Note that (5.1) and (5.2) are superseded, for $C$ large enough, by the bound
$$
\sum_{0<\g_n\le T;\g_{n+1}-\g_n>C/\log(T/2\pi)}1\;\ll\;
N(T)\exp(-AC)
\quad(A>0, C\ge C_0).\leqno(5.3)
$$
The bound (5.3) is Corollary 2 on p. 35 of A. Fujii \cite{5}. By (2.3) and
$C = \lambda \log(T/2\pi)$ with $\lambda >0$ sufficiently
large, the sum in (5.3) is empty. In that case
$$
N(T)\exp(-AC) = N(T)\exp\Bigl(-A\lambda \log(T/2\pi)\Bigr) = N(T)(T/2\pi)^{-A\lambda} < 1
$$
if $\lambda > 1/A, T\ge T_0>0$.
This shows that the bound in (5.3) is quite strong.

\medskip
%\break
\head
6. Sums of reciprocals of $\g_{n+1}-\g_n$.
\endhead

\medskip
The sum $S_k(T)$ in (1.1) clearly makes sense not only for $k>0$, but for $k<0$ as well (for $k=0$ the sum
is just $N(T)$, so it need not be considered). When $k<0$ one has obviously to assume the condition
$\g_{n+1} \ne \g_n$, or equivalently $\g_{n+1} > \g_n$, to avoid zeros in the denominator. Such a
condition is also natural when $k>0$, since
$$
(\g_{n+1}-\g_n)^k \equiv 0\qquad(k>0, \g_{n+1} = \g_n).
$$
There seem to be no results concerning $S_k(T)$ in the literature when $k<0$.   Even the sum
$S_{-1}(T)$ seems elusive.

\medskip
We shall consider here the somewhat less difficult sum
$$
H(T) := \sum_{0<t_n\le T,t_{n+1}\ne t_n}(t_{n+1}-t_n)^{-1},
$$
where $0 < t_1\le t_2\le t_3\le \ldots$ are the ordinates of zeta zeros on the critical line
$\Re s = \hf$, or equivalently, the zeros of Hardy's function $Z(t)$. Further let
$$
R(T) := \sum_{0<t_n\le T,t_{n+1}\ne t_n}1.
$$
If $\rho_n = \hf + it_n$ is a simple zero of $\z(s)$, then we cannot have $t_{n+1}=t_n$. Thus
$R(T)$ counts all simple zeros on the critical line, and the number of those for which $0<t_n\le T$
is $\gg T\log T$. In fact, H.M. Bui, B. Conrey and M.P. Young \cite{2} showed that more than 40\%
of the zeros counted by $N(T)$ are simple and on the critical line. More recently
N. Robles, A. Roy and  A. Zaharescu \cite{17} proved that at least 41.0725\% of the zeros of $\z(s)$
 are on the critical line and at least 40.5824\% of the zeros of $\z(s)$ are both on the critical line and
simple.

\medskip
Thus for some $D$ satisfying
$D>2/5$ we have
$$
R(T) \;\ge \frac {DT}{2\pi}\log \frac{T}{2\pi}\qquad(T\ge T_0 >0).\leqno(6.1)
$$
On the other hand, by using the Cauchy-Schwarz inequality, we obtain
$$
\eqalign{
R(T) &= \sum_{0<t_n\le T,t_{n+1}\ne t_n}\frac{1}{\sqrt{t_{n+1}-t_n}}\cdot \sqrt{t_{n+1}-t_n}
\cr&
\le \left\{H(T)\sum_{0<t_n\le T,t_{n+1}\ne t_n}(t_{n+1}-t_n)\right\}^{1/2}.\cr}\leqno(6.2)
$$
Since $t_{n+1}-t_n \ll t_n^{1/6}$ (see Chapter 9 of \cite{11}), it follows that
$$
\sum_{0<t_n\le T,t_{n+1}\ne t_n}(t_{n+1}-t_n) = T + o(T)\qquad(T\to\infty).\leqno(6.3)
$$
From (6.1)--(6.3) we obtain that
$$
\frac {DT}{2\pi}\log \frac{T}{2\pi} \le \sqrt{H(T)(T+o(T))},
$$
which gives

\medskip
THEOREM 3. {\it We have}
$$
H(T) = \sum_{0<t_n\le T,t_{n+1}\ne t_n}(t_{n+1}-t_n)^{-1} \ge
\frac{T}{(5\pi)^2}{\left(\log \frac{T}{2\pi}\right)}^2
\qquad(T\ge T_1>0).\leqno(6.4)
$$
An upper bound for $H(T)$ seems difficult to obtain.

\smallskip\smallskip\smallskip

We have
$$
\eqalign{
H(T)&\le \max_{0<t_n\le T,t_{n+1}\ne t_n}{(t_{n+1}-t_n)}^{-1}N(T)\cr&
\le \max_{0<t_n\le T,t_{n+1}\ne t_n}(t_{n+1}-t_n)^{-1}\left(\frac{T}{2\pi}\log \frac{T}{2\pi} + O(T)\right),
\cr}
$$
and from (6.4) it follows that we obtain
$$
\max_{0<t_n\le T,t_{n+1}\ne t_n}{(t_{n+1}-t_n)}^{-1} \ge \frac{2}{25\pi}\log \frac{T}{2\pi}
\qquad(T\ge T_1>0).\leqno(6.5)
$$
or equivalently
$$
\min_{0<t_n\le T,t_{n+1}\ne t_n}(t_{n+1}-t_n) \le \frac{25\pi}{2\log \frac{T}{2\pi}}
\qquad(T\ge T_1>0).\leqno(6.6)
$$
If one considers the analogous problem with the sequence $\{t_n\}$ replaced by the sequence
$\{\g_n\}$, then only the analogue of (6.1) is not obvious, namely
$$
\sum_{0<\g_n\le T,\g_{n+1}\ne \g_n}1 \;\gg\; T\log T. \leqno(6.7)
$$
However, the sum in (6.7) certainly counts simple zeros (with $\g_n\le T$)
on the critical line, and as already
mentioned, there are $\gg T\log T$ of these. Thus (6.7) holds,
and the  rest of the preceding argument easily carries through.
Alternatively, since the sum in (6.7) certainly also counts {\it distinct} zeros of $\z(s)$,
and there are at least 70\% of distinct zeta-zeros (see H. Ki and Y. Lee [15]), 
we can obtain an even better bound for the sum in (6.7).
Therefore we can obtain the analogues of (6.4)--(6.6)
for the sequence $\{\g_n\}$, with different explicit constants, of course.

\medskip
%\vfill
%\eject
%\topglue1cm
%\vskip2cm
%\medskip
\Refs
\medskip

\item{[1]} J. Bredberg, Large gaps between consecutive zeros, on the critical line,
of the Riemann zeta-function, preprint available at
\tt{arXiv:1101.3197v3}.

\smallskip\rm
\item{[2]} H.M. Bui, B. Conrey, and M.P. Young,
More than $41\%$ of the zeros of the zeta function are on the critical line,
 Acta Arith. {\bf150}(2011), no. 1, 35-64.

\smallskip
\item{[3]} E. Carneiro,  V. Chandee and
M.B. Milinovich, Bounding $S(t)$ and $S_1(t)$ on the Riemann hypothesis,
Math. Ann. {\bf356}(2013), 939-968.

\smallskip
\item{[4]} S. Feng and X. Wu, On gaps between zeros of the Riemann zeta-function,
J. Number Theory {\bf132}(2012), 1385-1397.

\smallskip
\item{[5]} A. Fujii, On the distribution of the zeros of the Riemann zeta-function in short
intervals, Bull. Amer. Math. Soc. {\bf 81}(1975), 139-142.

\smallskip
\item{[6]} A. Fujii, On the zeros of Dirichlet L-functions. II.
(With corrections to ``On the zeros of Dirichlet L-functions. I'', and the subsequent papers),
Trans. Am. Math. Soc. {\bf267}(1981), 33-40.

\smallskip
\item{[7]} A. Fujii, On the gaps between the consecutive zeros of the Riemann zeta-function,
Proc. Japan Acad. {\bf66}, Ser. A. (1990), 97-100.

\smallskip
\item{[8]}
D.A.  Goldston and S.M. Gonek,  A note on $S(t)$ and the zeros of the Riemann zeta-function,
 Bull. Lond. Math. Soc. {\bf39}(3)(2007),  482-486.

\smallskip
\item{[9]}
S.M. Gonek and A. Ivi\'c,  On the distribution of  positive and negative
values of Hardy's $Z$-function, to appear, preprint available at
\tt{arXiv:1604.00517.} \rm

\smallskip

\item{[10]} R.R. Hall and W.K. Hayman, Hyperbolic distance and distinct zeros of the
Riemann zeta-function in small regions, J. reine angew. Math. {\bf526}(2000), 35-59.

\smallskip
\item{[11]} A. Ivi\'c, The Riemann zeta-function, John Wiley \&
Sons, New York 1985 (reissue,  Dover, Mineola, New York, 2003).

\smallskip\rm
\item{[12]} A. Ivi\'c, On sums of gaps between the zeros of $\zeta(s)$ on the critical line,
      Univ. Beograd. Publ. Elektrotehn. Fak. Ser. Mat. {\bf6}(1995), 55-62.

\smallskip
\item{[13]} A. Ivi\'c, The theory of Hardy's $Z$-function, Cambridge University Press,
Cambridge, 2012, 245 pp.

\smallskip
\item{[14]} D. Joyner, On the Dyson--Montgomery hypothesis, Proc. Amalfi Conf. Analytic
Number Theory 1989, Univ. Salerno, 1992, 257-261.

\smallskip
\item{[15]} H. Ki and Y. Lee, Zeros of the derivatives of the Riemann zeta-function, Functiones et
Approximatio {\bf47}, No. 1(2012), 79-87.

\smallskip
\item{[16]} A.M. Odlyzko, On the Distribution of Spacings Between Zeros of the Zeta Function,
Math. Comp. Vol. {\bf48} No. {\bf177}(1987), 273-308.

\smallskip
\item{[17]} N. Robles, A. Roy and  A. Zaharescu,
Twisted second moments of the Riemann zeta-function and applications,
J. Math. Anal. Appl. {\bf434}(2016), no. 1, 271-314.

\smallskip
 \item{[l8]} A. Selberg, The zeta-function and the Riemann Hypothesis, Skandinaviske
 Mathematikerkongres {\bf10}(1946), 187-200.

\smallskip
\item{[19]} E.C. Titchmarsh, The theory of the Riemann zeta-function, 2nd ed. edited by
D.R. Heath-Brown, Oxford, Clarendon Press, 1986.

\smallskip
\item{[20]} T. Trudgian, An improved upper bound for the argument of the
Riemann zeta-function on the critical line II, J. Number Theory {\bf134}(2014), 280-292.
\vfill

\endRefs
%\bigskip

\enddocument

\bye